\documentclass[10pt, a4paper]{article}
\usepackage[pass]{geometry}

\usepackage{amssymb}
\usepackage{amsmath}
\usepackage{amsthm}
\usepackage{mathrsfs}
\usepackage{microtype}
\usepackage[utf8]{inputenc}
\usepackage[english]{babel}
\usepackage{url}
\usepackage{enumitem}
\setlist[enumerate]{label=\arabic*), nosep, leftmargin=*}

\makeatletter

\renewcommand{\maketitle}{%
    \begingroup
    \parindent=0pt
    \centering
    {\normalsize\bfseries \@author\par}
    \vspace{0.9em}
    {\large\bfseries \@title\par}
    \vspace{0.6em}
    {\small
    Moscow Institute of Physics and Technology,
    Dolgoprudny, Russia\par
    Email: \texttt{beketov.di@phystech.edu}\par}
    \vspace{1.2em}
    \endgroup
}

\renewenvironment{abstract}
{%
    \small
    \setlength{\emergencystretch}{2em}%
    \begin{list}{}{%
        \setlength{\leftmargin}{\parindent}%
        \setlength{\rightmargin}{\parindent}%
        \setlength{\topsep}{0pt}%
        \setlength{\partopsep}{0pt}%
        \setlength{\parsep}{0pt}%
        \setlength{\itemsep}{0pt}%
    }%
    \item[]\noindent\textbf{Abstract.}\ %
}
{%
    \end{list}%
}

\newenvironment{keywords}
{%
    \small
    \begin{list}{}{%
        \setlength{\leftmargin}{\parindent}%
        \setlength{\rightmargin}{\parindent}%
        \setlength{\topsep}{0pt}%
        \setlength{\partopsep}{0pt}%
        \setlength{\parsep}{0pt}%
        \setlength{\itemsep}{0pt}%
    }%
    \item[]\noindent\textbf{Keywords:}\ %
}
{%
    \end{list}%
}

\makeatother

\newtheoremstyle{compact}
  {6pt}{6pt}
  {\normalfont}{}
  {\bfseries}{.}
  {0.5em}{}
\newtheoremstyle{compactitalic}
  {6pt}{6pt}
  {\itshape}{}
  {\bfseries}{.}
  {0.5em}{}
\makeatletter
\renewenvironment{proof}[1][\proofname]{%
  \par\addvspace{6pt}\pushQED{\qed}\normalfont
  \noindent{\itshape #1\@addpunct{.}}\ \ignorespaces
}{%
  \popQED\par\addvspace{6pt}
}
\makeatother
\theoremstyle{compactitalic}
\newtheorem{theorem}{Theorem}
\newtheorem{lemma}[theorem]{Lemma}

\theoremstyle{compact}
\newtheorem{definition}[theorem]{Definition}
\newtheorem{example}[theorem]{Example}

\begin{document}

\setlength{\parindent}{1.5em}
\setlength{\parskip}{0pt}
\setlength{\emergencystretch}{2em}

\author{Dmitriy I. Beketov}
\title{Generalized Limits of Finitely Additive Measures on Locally Compact Groups}
\maketitle

\begin{abstract}
Let \(\mathscr{E}\) be an increasing family of Borel subsets of a locally compact topological group \(G\) and let \(\mu\) be a nonnegative finitely additive \(\sigma\)-finite measure on the Borel \(\sigma\)-algebra of \(G\) whose values on \(\mathscr{E}\) are finite but unbounded. Applying a generalized limit pointwise to the net of measures \(\left\{B\mapsto\frac{\mu(B\cap E)}{\mu(E)}\right\}_{E\in \mathscr{E}}\) yields a finitely additive measure \(\varphi_\mu\), which can be viewed as a generalized normalization of \(\mu\) to total mass one. To obtain a nontrivial necessary condition on \(\mu\) for the resulting measure \(\varphi_{\mu}\) to be invariant, we establish two results of independent interest: we extend the Yosida–Hewitt decomposition to finitely additive \(\sigma\)-finite measures and prove a theorem previously stated by Bingham and Ostaszewski concerning translates of sets of Haar measure zero. The first theorem yields a decomposition of \(\mu\) into a purely finitely additive component and two countably additive components, respectively singular and absolutely continuous with respect to Haar measure. Using the second result, we show that the invariance of \(\varphi_\mu\) implies that the singular component \(\mu_s\) is negligible relative to \(\mu\) and can be removed without changing \(\varphi_\mu\). A counterexample demonstrates that the purely finitely additive component \(\mu_p\) need not have these properties.
\end{abstract}

\begin{keywords}
finitely additive measures, invariant measures, Yosida--Hewitt decomposition, $\sigma$-finite measures, locally compact groups.
\end{keywords}

\section*{\centering Introduction}
Invariant means on locally compact groups (and the corresponding invariant finitely additive measures) find applications in operator theory \cite[Theorem 4.3]{GhaffariAmirjan2016}, additive combinatorics \cite[Theorems 1.4 and 1.7]{GriesmerLeLe2023}, quantum probability and statistics \cite[pp. 105--106]{Holevo2001}, \cite{Srinivas1980}, geometric measure theory \cite[Theorem 1.1]{ConnesMcDonaldSukochevZanin2019} and the study of dynamical systems \cite{Sakbaev2024}. One way to construct invariant finitely additive measures is to apply a generalized limit to a net of normalized restrictions of a measure to a suitable family of sets. In particular, this approach was used in \cite{AdlerHamilton1973} to prove the sufficiency of the F{\o}lner condition for amenability of a discrete semigroup. Much later, the same idea was applied to Lebesgue measure on the Borel $\sigma$-algebra of the real line \cite{Sakbaev2024}:
\begin{equation}\label{eq:Lebesgue-limit}
\nu(A) = \lim\limits_{x, \mathbb{F}} \frac{\lambda(A \cap [-x, x])}{2x},
\end{equation}
where \( \mathbb{F} \) is an ultrafilter on \( (0; +\infty) \) containing no bounded sets \cite[Ch. 1, \S\S\ 6 and 7]{Bourbaki1998}. The resulting measure was used to obtain a unitary representation of the group of infinite-dimensional Hamiltonian flows \cite{Sakbaev2024} and was studied in greater detail in \cite{DzhenzherSakbaev2026}. In the present paper, we generalize construction \eqref{eq:Lebesgue-limit} to arbitrary locally compact topological groups and finitely additive $\sigma$-finite measures. Given an increasing family of Borel sets \( \mathscr{E} = \{ E_{\alpha}\}_{\alpha \in I} \) and a generalized limit \( \Lambda \) on \(I\), we consider the normalized measure:
\[
\varphi_{\mu}(A) = \Lambda \frac{\mu(A \cap E_{\alpha})}{\mu(E_{\alpha})}.
\]\par
A central aim of this paper is to obtain a nontrivial necessary condition on \(\mu\) for the resulting measure \(\varphi_\mu\) to be invariant. To this end, we first show that \(\mu\) admits a decomposition into three \(\sigma\)-finite measures: a countably additive measure \(\mu_{ac}\) absolutely continuous with respect to Haar measure, a countably additive measure \(\mu_s\) singular with respect to Haar measure, and a purely finitely additive measure \(\mu_p\). We then prove that, if \(\varphi_\mu\) is invariant, the singular component \(\mu_s\) is negligible relative to \(\mu\) in the sense that \( \Lambda\frac{\mu_s(E_\alpha)}{\mu(E_\alpha)}=0 \). From this condition, we further derive the identity \(\varphi_\mu=\varphi_{\mu-\mu_s}\). We use these two results to show, by means of a counterexample, that neither conclusion extends in general to the purely finitely additive component \( \mu_p \).\par
The body of the paper is divided into three sections. In Section \ref{sec:singular-measures}, we prove a theorem on measures on locally compact topological groups that was stated in \cite{BinghamOstaszewski2020} but proved there only under an additional assumption. This theorem will be used later to study singular measures. Then, in Section 2, we generalize the Yosida–Hewitt decomposition to finitely additive \(\sigma\)-finite measures that may take the value \(+\infty\), using the known existence and uniqueness theorem for the decomposition of finite-valued finitely additive measures (see \cite{Diestel1970}; a proof can be found in \cite[Theorem 8.16]{Konig1997}). The above decomposition of \( \mu \) into three components will follow from this result and the Lebesgue decomposition theorem \cite[Theorem 15.14]{Nielsen1997}. Finally, in Section \ref{sec:generalized-limits}, we use the results of the first two sections to establish the assertions concerning \(\mu_s\) and \(\mu_p\) outlined in the preceding paragraph.\par

\section{\centering Singular Measures on Topological Groups}\label{sec:singular-measures}
In this paper, by a measure we mean a finitely additive set function defined on a ring \( \mathcal{R} \) of subsets of a set \( S \) and taking values in \( [0; +\infty] \).\par
Let \( X \) be a topological space. We use the following notation:
\begin{enumerate}
\item \( \chi_E \) is the indicator function of a set \( E \subset X \);
\item \( \mathscr{B}(X) \) is the Borel $\sigma$-algebra on \( X \);
\item \( Z(X) = \{ f^{-1}(\{0\}): f \text{ is continuous on } X \} \) is the family of all zero sets of \( X \);
\item \( \mathcal{B}a(X) \) is the Baire $\sigma$-algebra, i.e., the $\sigma$-algebra generated by \( Z(X) \).
\end{enumerate}\par
The elements of the Baire $\sigma$-algebra will be called Baire sets. A countably additive measure \( \mu: \mathscr{B}(X) \to [0; +\infty] \) is called completion regular if, for every \( A \in \mathscr{B}(X) \), we have \( \mu(A) = \sup\{ \mu(B): B \in Z(X), B \subset A \} \).\par
Fix a locally compact topological group \( G \). Denote its left Haar measure by \( \lambda \). Define the modular function by the identity \( \lambda(Ex) = \Delta(x) \lambda(E) \) for all \( x \in G \) and \( E \in \mathscr{B}(G) \). It is well known (see, e.g., \cite[Theorem 9.11.6]{Bogachev2007}) that left Haar measure is completion regular.
\begin{theorem}\label{thm:SakM}
\textnormal{\cite[Theorem SakM]{BinghamOstaszewski2020}} Let \( \nu \) be a countably additive $\sigma$-finite measure on \( \mathscr{B}(G) \), and let \( S \) belong to \( \mathscr{B}(G) \) and satisfy \( \lambda(S)=0 \). Then, for \( \lambda \)-almost every \( t \in G \), we have \( \nu(tS)=0 \).
\end{theorem}
In \cite{BinghamOstaszewski2020}, this theorem is deduced from a result proved only for Polish groups. We therefore give our own proof for the general case, based on a similar approach.\par
\begin{proof}[Proof]\par\nobreak
1. First suppose that \( G \) is $\sigma$-compact. Completion regularity and $\sigma$-finiteness of \( \lambda \) imply that it is outer regular on the Borel $\sigma$-algebra with respect to the Baire $\sigma$-algebra. In particular, there exists \( N \in \mathcal{B}a(G) \) such that \( S \subset N \) and \( \lambda(N)=0 \).\par
\begin{samepage}
By \cite[Lemma 6.3.6]{Bogachev2007}, for arbitrary topological spaces \( X \) and \( Y \), a continuous map \( f: X \to Y \), and a Baire set \( B \subset Y \), the inverse image \( f^{-1}(B) \) is also a Baire set. Taking \( X = G \times G \), \( Y = G \), \( f: (t, x) \mapsto t^{-1}x \), and \( B=N \), we conclude that \( E = \{ (t, x): t^{-1}x \in N \} \) is a Baire set.\par
\end{samepage}
Choose a family of compact sets \( \{K_n\}_{n=1}^{\infty} \) such that \( G = \bigcup\limits_{n=1}^{\infty} K_n \). Then \( G \times G = \bigcup\limits_{n,m \in \mathbb{N}} (K_n \times K_m) \). For any \( n, m \in \mathbb{N} \), by the Tychonoff theorem, \( K_n \times K_m \) is compact in \( G \times G \). Thus \( G \times G \) is $\sigma$-compact. Consequently, both \( G \) and \( G \times G \) are Lindel\"of spaces. Hence \( \mathcal{B}a(G \times G) = \mathcal{B}a(G) \otimes \mathcal{B}a(G) \) \cite[Proposition 6.10.7]{Bogachev2007}, so \( E \in \mathcal{B}a(G) \otimes \mathcal{B}a(G) \). Since the Baire $\sigma$-algebra is contained in the Borel $\sigma$-algebra, \( E \in \mathscr{B}(G) \otimes \mathscr{B}(G) \).\par
By Fubini's theorem \cite[Theorem 8.8]{Rudin1987},
\[
\begin{gathered}
\int\limits_G \nu(tN) \, d\lambda(t) =
\int\limits_G \left( \int\limits_G \chi_E(t, x)
\, d\nu(x) \right) \, d\lambda(t)\\
=\int\limits_G \left( \int\limits_G \chi_E(t, x)
\, d\lambda(t) \right) \, d\nu(x) =
\int\limits_G \lambda(xN^{-1}) \, d\nu(x).
\end{gathered}
\]
Using the equalities
\[
\lambda(xN^{-1}) = \lambda(N^{-1}) = \int\limits_N \Delta(g^{-1})\, d\lambda(g) = 0
\]
\cite[Proposition 2.31]{Folland2016}, we obtain \( \int\limits_G \nu(tN)\, d\lambda(t) = \int\limits_G 0\, d \nu(x) = 0 \). Hence \( \nu(tN) = 0 \) for \( \lambda \)-almost every \( t \in G \). Since \( \nu(tS) \leq \nu(tN) \), it follows that \( \nu(tS)=0 \) for \( \lambda \)-almost every \( t \in G \).\par
Note that the Borel measurability of \( S \) and the continuity of \( f \) imply that \( E' = \{ (t, x): t^{-1}x \in S \} \) belongs to \( \mathscr{B}(G \times G) \). However, the $\sigma$-algebras \( \mathscr{B}(G \times G) \) and \( \mathscr{B}(G) \otimes \mathscr{B}(G) \) do not coincide in general. Thus Fubini's theorem cannot be applied directly to the function \( \chi_{E'} \) without considering the Baire $\sigma$-algebras.\par
2. Now consider an arbitrary locally compact group \( G \). By \cite[Proposition 2.4]{Folland2016}, \( G \) has an open, closed, $\sigma$-compact subgroup \( G_0 \). Choose \( Y \subset G \) such that \( G = \bigsqcup\limits_{y \in Y} (yG_0) \). Then the set \( A = \{a \in Y: S \cap (aG_0) \neq \emptyset \} \) is at most countable \cite[Proposition 2.22]{Folland2016}. Moreover, $\sigma$-finiteness of \( \nu \) implies that \( B = \{ b \in Y: \nu(bG_0) > 0 \} \) is also at most countable.\par
Suppose that for some element \( t \) of \( G \) we have \( \nu(tS) > 0 \). By countable additivity of \( \nu \), there exists \( a_t \in A \) such that \( \nu(t(S \cap a_tG_0)) > 0 \). Since \( \nu(ta_tG_0) \geq \nu(t(S \cap a_tG_0)) > 0 \), there exists \( b_t \in B \) such that \( ta_tG_0 = b_tG_0 \). Consequently, \( ta_t \in b_tG_0 \), or equivalently, \( t \in b_tG_0a_t^{-1} \). Thus
\begin{equation}\label{eq:reduction-to-subgroup}
\begin{gathered}
t \in \{ g \in b_tG_0a_t^{-1}: \nu(g(S \cap a_tG_0)) > 0 \}\\
=\{ b_tha_t^{-1}: h \in G_0, \nu(b_tha_t^{-1}(S \cap a_t G_0)) > 0 \}\\
=b_t\{ h \in G_0: \nu(b_tha_t^{-1}(S \cap a_tG_0)) > 0 \}a_t^{-1}.
\end{gathered}
\end{equation}\par
For arbitrary \( a \in A \) and \( b \in B \), define the set
\[
N_{a,b} = \{ h \in G_0: \nu(bha^{-1}(S \cap aG_0)) > 0 \}
\]
and the measure \( \nu_b \) by setting, for every \( E \in \mathscr{B}(G_0) \), \( \nu_b(E) = \nu(bE) \). Applying the $\sigma$-compact case already proved to the topological group \( G_0 \), the measure \( \nu_b \), and the set \( a^{-1}(S \cap aG_0) \), we obtain \( \lambda(M_{a, b}) = 0 \) for some \( M_{a,b} \in \mathscr{B}(G_0) \) containing \( N_{a,b} \). Then \( \lambda(b M_{a, b} a^{-1}) = \lambda(M_{a,b} a^{-1}) = \Delta(a^{-1}) \lambda(M_{a,b}) = 0 \). By (\ref{eq:reduction-to-subgroup}), we have the inclusions
\[
\{ t \in G: \nu(tS) > 0 \} \subset \bigcup\limits_{a \in A} \bigcup\limits_{b \in B} ( bN_{a,b}a^{-1}) \subset \bigcup\limits_{a \in A} \bigcup\limits_{b \in B} ( bM_{a,b}a^{-1}).
\]
Therefore, the inequalities
\[
\lambda \left( \bigcup\limits_{a \in A} \bigcup\limits_{b \in B} ( bM_{a,b}a^{-1} ) \right) \leq \sum\limits_{a \in A} \sum\limits_{b \in B} \lambda(bM_{a,b}a^{-1}) = 0
\]
prove the assertion.
\end{proof}
\begin{lemma}\label{lem:disjoint-translates}
Let \( \gamma \in ca_{\sigma}^+(G, \mathscr{B}(G)) \), and suppose that \( \gamma \) and \( \lambda \) are singular. Then, for every \( m \in \mathbb{N} \), there exist a set \( H \in \mathscr{B}(G) \) and group elements \( g_1, ..., g_m \in G \) such that \( \lambda(H)=0 \), \( \gamma \left( G \setminus H \right) = 0 \), and the sets \( g_1H, ..., g_mH \) are pairwise disjoint.
\end{lemma}
\begin{proof}[Proof]\par\nobreak
1. The singularity of \( \gamma \) and \( \lambda \) means that there exists \( A \in \mathscr{B}(G) \) such that \( \lambda(A)=0 \) and \( \gamma \left( G \setminus A \right) = 0 \). Denote the set of all translations \( g \in G \) satisfying \( \gamma(gA)=0 \) by \( \mathcal{G} \). By Theorem \ref{thm:SakM}, \( \lambda \left( G \setminus \mathcal{G} \right) = 0 \).\par
2. We prove by induction that for every \( m \in \mathbb{N} \) there exist group elements \( g_1, ..., g_m \in G \) such that for all \( i, j \leq m \) with \( i < j \) we have \( g_i^{-1} g_j \in \mathcal{G} \). The base case is immediate: \( g_1 \in G \) can be chosen arbitrarily. Assume that the assertion holds for \( m=k-1 \), where \( k \geq 2 \). We prove it for \( m=k \). Consider elements \( g_1, ..., g_{k-1} \in G \) whose existence is guaranteed by the induction hypothesis.\par
Define \( B = \bigcap\limits_{i=1}^{k-1} (g_i \mathcal{G}) \). Since
\[
\lambda \left( G \setminus B \right) \leq \sum\limits_{i=1}^{k-1} \lambda \left( G \setminus (g_i\mathcal{G}) \right) = \sum\limits_{i=1}^{k-1} \lambda \left( G \setminus \mathcal{G} \right) = 0,
\]
the set \( B \) is nonempty. Choose any \( g_k \in B \). Clearly, for every \( i < k \) we have \( g_i^{-1} g_k \in \mathcal{G} \).\par
3. For any \( m \in \mathbb{N} \), consider the elements \( g_1, ..., g_m \in G \) obtained in the preceding step. We show that
\[
H = A \setminus \left( \bigcup\limits_{1 \leq i < j \leq m} (g_i^{-1} g_j A) \right)
\]
is the required set.\par
From \( H \subset A \) and \( \lambda(A)=0 \), it follows that \( \lambda(H) = 0 \). Moreover,
\begin{gather*}
\gamma(G\setminus H)
=\gamma(G\setminus A)+\gamma\left(A\cap\left(\bigcup_{1\leq i<j\leq m}(g_i^{-1}g_j A)\right)\right)\\
\leq\gamma\left(\bigcup_{1\leq i<j\leq m}(g_i^{-1}g_j A)\right)
\leq\sum_{1\leq i<j\leq m}\gamma(g_i^{-1}g_j A)=0.
\end{gather*}
Now take any \( i, j \in \{1, 2, ..., m\} \) with \( i < j \) and suppose that \( g_iH \cap g_jH \neq \emptyset \). Then there exist \( x, y \in H \) such that \( g_i x = g_j y \). Consequently, \( x = g_i^{-1} g_j y \). Hence \( H \cap (g_i^{-1} g_j A) \supset \{ x \} \cap \{ g_i^{-1}g_jy \} = \{ x \} \). But \( H \cap (g_i^{-1} g_j A) = \emptyset \). This contradiction implies that \( g_i H \cap g_j H = \emptyset \) whenever \( i \neq j \).
\end{proof}

\section{\centering The Yosida--Hewitt Decomposition of a $\sigma$-Finite Measure}\label{sec:Yosida-Hewitt}
Fix a ring \( \mathcal{R} \) of subsets of a set \( S \). We use the following notation:
\begin{enumerate}
\item \( ba^+(S, \mathcal{R}) \) is the set of all bounded measures on \( \mathcal{R} \);
\item \( ca^+(S, \mathcal{R}) \) is the set of all bounded countably additive measures on \( \mathcal{R} \);
\item For \( B \in \mathcal{R} \), \( \mu \llcorner B \) denotes the measure on \( \mathcal{R} \) defined, for every \( A \in \mathcal{R} \), by \( (\mu \llcorner B)(A) = \mu(A \cap B) \).
\end{enumerate}\par
If \( \mathcal{R} \) is an algebra and \( \mu \colon \mathcal{R} \to \mathbb{R} \) is a bounded finitely additive set function, we define the total variation of \( \mu \) as follows:
\[
\| \mu \| = \sup\limits_{n \in \mathbb{N}, \{ E_i\}_{i=1}^n \subset \mathcal{R}: \bigsqcup\limits_{i=1}^n E_i = S } \sum\limits_{i=1}^n \left|\mu(E_i)\right|.
\]
\begin{definition}\label{def:sigma-finite}
A measure \( \mu \) on \( \mathcal{R} \) is called $\sigma$-finite if there exists a countable family \( \{B_j\}_{j=1}^{\infty} \) of elements of \( \mathcal{R} \) such that \( S = \bigsqcup\limits_{j=1}^{\infty} B_j \) and, for every \( j \in \mathbb{N} \), we have \( \mu(B_j) < +\infty \). We denote the set of all $\sigma$-finite measures on \( \mathcal{R} \) by \( ba_{\sigma}^+(S, \mathcal{R}) \), and the set of all countably additive $\sigma$-finite measures on \( \mathcal{R} \) by \( ca_{\sigma}^+(S, \mathcal{R}) \).
\end{definition}
\begin{definition}\label{def:purely-finitely-additive}
A measure \( \mu \in ba_{\sigma}^+(S, \mathcal{R}) \) is called purely finitely additive if every countably additive measure \( \gamma \) on \( \mathcal{R} \) satisfying \( \gamma \leq \mu \) is such that \( \gamma = 0 \). We denote the set of purely finitely additive measures on \( \mathcal{R} \) by \( \mathrm{pfa}_{\sigma}^+(S, \mathcal{R}) \), and the set of bounded purely finitely additive measures on \( \mathcal{R} \) by \( \mathrm{pfa}^+(S, \mathcal{R}) \).
\end{definition}
Every measure \( \mu \in ba_{\sigma}^+(S, \mathcal{R}) \) that does not take the value \( +\infty \) admits a unique decomposition into measures \( \mu_c \in ca_{\sigma}^+(S, \mathcal{R}) \) and \( \mu_p \in \mathrm{pfa}_{\sigma}^+(S, \mathcal{R}) \) that likewise do not take the value \( +\infty \) \cite[Theorem 8.16]{Konig1997}. This decomposition of \( \mu \) is called the Yosida--Hewitt decomposition.
\begin{theorem}\label{thm:Yosida-Hewitt}
Every \( \mu \in ba_{\sigma}^+(S, \mathcal{R}) \) admits a decomposition into measures \( \mu_c \in ca_{\sigma}^+(S, \mathcal{R}) \) and \( \mu_p \in \mathrm{pfa}_{\sigma}^+(S, \mathcal{R}) \), uniquely characterized by the following property: for every pair \( \tilde\mu_c \in ca_{\sigma}^+(S, \mathcal{R}) \), \( \tilde\mu_p \in \mathrm{pfa}_{\sigma}^+(S, \mathcal{R}) \) satisfying \( \mu = \tilde\mu_c + \tilde\mu_p \), we have \( \tilde\mu_c = \mu_c \) and \( \mu_p \leq \tilde\mu_p \).
\end{theorem}
\begin{proof}[Proof]\par\nobreak
1. Let \( r \) denote the ring of sets on which \( \mu \) takes finite values: \( r = \{ A \in \mathcal{R}: \mu(A) < +\infty \} \). Consider the Yosida--Hewitt decomposition of \( \mu \) on the ring \( r \): \( \mu |_r = \mu_c^r + \mu_p^r \).\par
2. For \( \mu \), take a family \( \{B_j\}_{j=1}^{\infty} \) from the definition of $\sigma$-finiteness. Note that every \( A \) in \( \mathcal{R} \) is a countable union of elements of \( r \): \( A = \bigcup\limits_{j=1}^{\infty} (A \cap B_j) \). Thus the ring \( \mathcal{R} \) is contained in the $\sigma$-ring generated by \( r \). Hence, by the Carath\'eodory extension theorem \cite[p. 54]{Halmos1950}, \( \mu_c^r \) has a unique countably additive extension \( \mu_c \) to \( \mathcal{R} \).\par
3. Let \( \mathcal{R}_f \) be the ring of sets on which \( \mu_c \) is finite: \( \mathcal{R}_f = \{ A \in \mathcal{R}: \mu_c(A) < +\infty \} \). Note that for every \( A \in \mathcal{R}_f \) we have
\[
\mu_c(A) = \sum\limits_{j=1}^{\infty} \mu_c(A \cap B_j) = \sum\limits_{j=1}^{\infty} \mu_c^r(A \cap B_j) \leq \sum\limits_{j=1}^{\infty} \mu(A \cap B_j) \leq \mu(A).
\]
Consequently, on \( \mathcal{R}_f \) the measure \( \mu_p^f = \mu |_{\mathcal{R}_f} - \mu_c |_{\mathcal{R}_f} \) is well defined and extends \( \mu_p^r \).\par
4. On \( \mathcal{R} \), define a set function \( \mu_p \) as follows:
\[
\mu_p(E) = \sup \left\{ \mu_p^f(A): A \in \mathcal{R}_f, A \subset E \right\}.
\]
For any disjoint sets \( E_1, E_2 \in \mathcal{R} \), we have
\begin{align*}
&\mathrel{\phantom{=}}\mu_p(E_1 \sqcup E_2)\\
&=\sup \left\{ \mu_p^f(A): A \in \mathcal{R}_f, A \subset E_1 \sqcup E_2 \right\}\\ &=\sup\limits \left\{ \mu_p^f(A_1) + \mu_p^f(A_2): A_1, A_2 \in \mathcal{R}_f, A_1 \subset E_1, A_2 \subset E_2 \right\}\\
&=\sup \left\{ \mu_p^f(A_1): A_1 \in \mathcal{R}_f, A_1 \subset E_1 \right\} + \sup \left\{ \mu_p^f(A_2): A_2 \in \mathcal{R}_f, A_2 \subset E_2 \right\}\\
&=\mu_p(E_1)+\mu_p(E_2).
\end{align*}
Thus \( \mu_p \) is a measure on \( \mathcal{R} \).\par
5. Since \( \mu_p \) extends \( \mu_p^f \), for \( A \in \mathcal{R}_f \) we have \( \mu(A) = \mu_c(A) + \mu_p(A) \). If \( A \in \mathcal{R} \setminus \mathcal{R}_f \), then \( \mu(A) = +\infty = \mu_c(A) = \mu_c(A) + \mu_p(A) \). Therefore \( \mu = \mu_c + \mu_p \). If a countably additive measure \( \gamma \) on \( \mathcal{R} \) satisfies \( \gamma \leq \mu_p \), then \( \gamma |_r \leq \mu_p |_r = \mu_p^r \). Since \( \mu_p^r \) is purely finitely additive, \( \gamma |_r = 0 \). By the Carath\'eodory extension theorem, \( \gamma = 0 \). Hence \( \mu_p \) is purely finitely additive.\par
6. Consider arbitrary measures \( \tilde\mu_c \in ca_{\sigma}^+(S, \mathcal{R}) \) and \( \tilde\mu_p \in \mathrm{pfa}_{\sigma}^+(S, \mathcal{R}) \) satisfying \( \mu = \tilde\mu_c+\tilde\mu_p \). We show that \( \tilde\mu_p |_r \) is purely finitely additive. Indeed, suppose that a measure \( \gamma^r \in ca_{\sigma}^+(S, r) \) satisfies \( \gamma^r \leq \tilde\mu_p |_r \). By the Carath\'eodory extension theorem, \( \gamma^r \) admits a countably additive extension \( \gamma \) to \( \mathcal{R} \), and for every \( A \in \mathcal{R} \) we have
\[
\gamma(A) = \sum\limits_{j=1}^{\infty} \gamma^r(A \cap B_j) \leq \sum\limits_{j=1}^{\infty} \tilde\mu_p(A \cap B_j) \leq \tilde\mu_p(A),
\]
that is, \( \gamma \leq \tilde\mu_p \). Then \( \gamma=0 \), so \( \gamma^r=0 \). Therefore \( \tilde\mu_p |_r \in \mathrm{pfa}_{\sigma}^+(S, r) \).\par
7. By uniqueness of the Yosida--Hewitt decomposition for real-valued measures, we have \( \tilde\mu_c |_r = \mu_c |_r \) and \( \tilde\mu_p |_r = \mu_p |_r \). Again by the Carath\'eodory extension theorem, \( \tilde\mu_c = \mu_c \). Consequently, \( \tilde\mu_p |_{\mathcal{R}_f} = \mu |_{\mathcal{R}_f} - \tilde\mu_c |_{\mathcal{R}_f} = \mu |_{\mathcal{R}_f} - \mu_c |_{\mathcal{R}_f} = \mu_p |_{\mathcal{R}_f} \). Thus, for every \( E \in \mathcal{R} \), we have
\begin{gather*}
\mu_p(E) = \sup \left\{ \mu_p(A): A \in \mathcal{R}_f, A \subset E \right\}\\
= \sup \left\{ \tilde\mu_p(A): A \in \mathcal{R}_f, A \subset E \right\} \leq \tilde\mu_p(E).
\end{gather*}
Hence \( \mu_p \leq \tilde\mu_p \).
\end{proof}
\begin{example}\label{ex:nonunique-decomposition}
Define \( \mu_p \in \mathrm{pfa}_{\sigma}^+(\mathbb{R}, \mathscr{B}(\mathbb{R})) \) as follows: if \( E \in \mathscr{B}(\mathbb{R}) \) has finite Lebesgue measure, set \( \mu_p(E)=0 \); otherwise set \( \mu_p(E)=+\infty \). Note that Lebesgue measure on the algebra \( \mathscr{B}(\mathbb{R}) \) admits two distinct decompositions into measures \( \mu_c \in ca_{\sigma}^+(\mathbb{R}, \mathscr{B}(\mathbb{R})) \) and \( \mu_p \in \mathrm{pfa}_{\sigma}^+(\mathbb{R}, \mathscr{B}(\mathbb{R})) \): \( \lambda = \lambda + \mu_p = \lambda + 0 \). Thus the decomposition of a measure \( \mu \in ba_{\sigma}^+(S, \mathcal{R}) \) into \( \mu_c \in ca_{\sigma}^+(S, \mathcal{R}) \) and \( \mu_p \in \mathrm{pfa}_{\sigma}^+(S, \mathcal{R}) \) is not unique in general.
\end{example}
\begin{definition}\label{def:Yosida-Hewitt}
We call the decomposition of a measure \( \mu \in ba_{\sigma}^+(S, \mathcal{R}) \) whose existence and uniqueness are established in Theorem \ref{thm:Yosida-Hewitt} the Yosida--Hewitt decomposition.
\end{definition}

\section{\centering Generalized Limits of Measures}\label{sec:generalized-limits}
Let \( G \) be a locally compact topological group. Fix a linearly ordered set \( I \) with no greatest element, and a family of sets \( \mathscr{E} = \{ E_{\alpha} \}_{\alpha \in I} \subset \mathscr{B}(G) \) satisfying \( E_{\alpha} \subset E_{\beta} \) whenever \( \alpha < \beta \). Denote by \( \mathcal{E} \) the union of all sets in \( \mathscr{E} \), and by \( B(I) \) the vector space of bounded functions from \( I \) to \( \mathbb{R} \), equipped with the supremum norm.
\begin{definition}\label{def:positive-functional}
\cite[p. 34]{Rudin1987} A functional \( \Lambda \in B^*(I) \) is called positive if, for every nonnegative function \( f \in B(I) \), we have \( \Lambda f \geq 0 \).
\end{definition}
Consider the filter \( \mathfrak{F} \) generated by the base \( \mathfrak{G} = \{ \{ \alpha \in I: \alpha > \alpha_0 \}: \alpha_0 \in I \} \), and an ultrafilter \( \mathbb{F} \) extending \( \mathfrak{F} \) (see, e.g., \cite[Ch. 1, \S\S\ 6 and 7]{Bourbaki1998}). Let \( L \) be the set of functions \( f \in B(I) \) for which \( \lim\limits_{\mathfrak{F}} f \) exists. Extend \( \lim\limits_{\mathfrak{F}} \) to a positive functional \( \Lambda \in B^*(I) \). Such an extension can be given, for example, by \( \lim\limits_{\mathbb{F}} \).\par
Let \( J \subset I \) be such that for some \( \beta \in I \) we have \( \{ \alpha \in I: \alpha > \beta \} \subset J \). Also consider a bounded function \( f: J \to \mathbb{R} \). Suppose \( f_1 \) and \( f_2 \) are bounded extensions of \( f \) to \( I \). Then \( \Lambda f_1 - \Lambda f_2 = \Lambda(f_1 - f_2) = \lim\limits_{\mathfrak{F}}(f_1 - f_2) = 0 \). Thus the value of \( \Lambda g \) does not depend on the choice of a bounded extension \( g \) of \( f \) to \( I \). Denote the set of all such functions \( f \) by \( \tilde{B}(I) \). For \( f \in \tilde{B}(I) \), we use \( \Lambda f \) to mean the value of \( \Lambda \) on any bounded extension of \( f \) to \( I \).\par
Throughout the rest of the paper, all measures under consideration are assumed to be defined on the $\sigma$-algebra \( \mathscr{B}(G) \).\par
\begin{definition}\label{def:invariant-measure}
For any measure \( \mu \) and any element \( g \) of the group \( G \), define the measure \( L_g \mu \) by \( (L_g \mu)(A) = \mu(gA) \). The measure \( \mu \) is called invariant if for every \( g \in G \) we have \( L_g \mu = \mu \).
\end{definition}
\begin{definition}\label{def:admissible-measure}
A $\sigma$-finite measure \( \mu \) is called admissible if
\begin{enumerate}
\item for every \( \alpha \in I \) we have \( \mu(E_{\alpha}) < +\infty \);
\item \( \lim\limits_{\mathfrak{F}} \mu(E_{\alpha}) = +\infty \).
\end{enumerate}
\end{definition}
Condition 2) of this definition implies that, for every admissible measure \( \mu \) and every \( A \in \mathscr{B}(G) \), the map \( \alpha \mapsto \frac{\mu(A \cap E_{\alpha})}{\mu(E_{\alpha})} \) belongs to \( \tilde{B}(I) \). This allows us to give the following definition.\par
\begin{definition}\label{def:generalized-limit-measure}
For every admissible measure \( \mu \), define the measure \( \varphi_{\mu} \) by
\[
\varphi_{\mu}(A) = \Lambda \frac{\mu(A \cap E_{\alpha})}{\mu(E_{\alpha})}
\]
for all \( A \in \mathscr{B}(G) \).
\end{definition}
\begin{example}\label{ex:counting-measure}
Let \( G = \langle \mathbb{Z}, +\rangle \) have the discrete topology, let \( I=\mathbb{N} \) and \( \mathscr{E} = \{ \{-n, -n+1, ..., n-1, n \} \}_{n \in \mathbb{N}} \), let \( \mu \) be counting measure on \( \mathbb{Z} \), and let \( \Lambda = B \) be any Banach limit (see, e.g., \cite[III, \S\ 7]{Conway1985}). Then
\[
\varphi_{\mu}(A) = B \left( \frac{1}{2n+1} |\{a \in A: |a| \leq n \}| \right),
\]
where \( | \cdot | \) denotes cardinality.
\end{example}
\begin{example}\label{ex:Lebesgue-measure}
Let \( G = \langle \mathbb{R}, + \rangle \) have its usual topology, let \( I=(0;+\infty) \), \( \mathscr{E} = \{ [-x,x] \}_{x>0} \), \( \mu = \lambda \), and \( \Lambda = \lim\limits_{\mathbb{F}} \). Then
\[
\varphi_{\mu}(A) = \lim\limits_{x, \mathbb{F}} \frac{\lambda(A \cap [-x,x])}{2x}.
\]
As noted above, measures of this form were used in \cite{Sakbaev2024} and studied in \cite{DzhenzherSakbaev2026}.
\end{example}
\begin{definition}\label{def:equivalent-measures}
Two admissible measures \( \mu_1\) and \( \mu_2 \) are called equivalent if \( \varphi_{\mu_1} = \varphi_{\mu_2} \). We denote equivalence of measures by \( \mu_1 \sim \mu_2 \).
\end{definition}
\begin{lemma}\label{lem:sufficient-for-equivalence}
If admissible measures \( \mu \) and \( \gamma \) satisfy
\[
\Lambda \frac{ \left\| \mu \llcorner E_{\alpha} - \gamma \llcorner E_{\alpha} \right\| }{\mu(E_{\alpha})} = 0,
\]
then \( \mu \sim \gamma \).
\end{lemma}
\begin{proof}[Proof]\par\nobreak
For every \( \alpha \in I \) such that \( \mu(E_{\alpha}) > 0 \) and \( \gamma(E_{\alpha}) > 0 \), we have
\begin{gather*}
\left\| \frac{\mu \llcorner E_{\alpha}}{\mu(E_{\alpha})} - \frac{\gamma \llcorner E_{\alpha}}{\gamma(E_{\alpha})} \right\| \leq \left\| \frac{\mu \llcorner E_{\alpha}}{\mu(E_{\alpha})} - \frac{\gamma \llcorner E_{\alpha}}{\mu(E_{\alpha})} \right\| + \left\| \frac{\gamma \llcorner E_{\alpha}}{\mu(E_{\alpha})} - \frac{\gamma \llcorner E_{\alpha}}{\gamma(E_{\alpha})} \right\|\\
=\frac{\left\| \mu \llcorner E_{\alpha} - \gamma \llcorner E_{\alpha} \right\|}{\mu(E_{\alpha})} + \gamma(E_{\alpha}) \frac{|\mu(E_{\alpha})-\gamma(E_{\alpha})|}{\mu(E_{\alpha}) \gamma(E_{\alpha})} \leq 2 \frac{\left\| \mu \llcorner E_{\alpha} - \gamma \llcorner E_{\alpha} \right\|}{\mu(E_{\alpha})}.
\end{gather*}
Consequently, for every set \( A \in \mathscr{B}(G) \) we have
\begin{gather*}
\left| \varphi_{\mu}(A) - \varphi_{\gamma}(A) \right| = \left| \Lambda \left( \frac{\mu(A \cap E_{\alpha})}{\mu(E_{\alpha})} - \frac{\gamma(A \cap E_{\alpha})}{\gamma(E_{\alpha})} \right) \right|\\
\leq \Lambda \left| \frac{\mu(A \cap E_{\alpha})}{\mu(E_{\alpha})} - \frac{\gamma(A \cap E_{\alpha})}{\gamma(E_{\alpha})} \right| \leq \Lambda \left\| \frac{\mu \llcorner E_{\alpha}}{\mu(E_{\alpha})} - \frac{\gamma \llcorner E_{\alpha}}{\gamma(E_{\alpha})} \right\|\\
\leq 2 \Lambda \frac{ \left\| \mu \llcorner E_{\alpha} - \gamma \llcorner E_{\alpha} \right\| }{\mu(E_{\alpha})} = 0. \qedhere
\end{gather*}
\end{proof}
Consider the Yosida--Hewitt decomposition of an arbitrary admissible measure \( \mu \in ba_{\sigma}^+(G, \mathscr{B}(G)) \): \( \mu = \mu_p + \mu_c \). In turn, decompose \( \mu_c \) into the sum of a measure \( \mu_{ac} \) absolutely continuous with respect to left Haar measure and a measure \( \mu_s \) singular with respect to left Haar measure \cite[Theorem 15.14]{Nielsen1997}.
\begin{theorem}\label{thm:remove-singular-part}
If \( \varphi_{\mu} \) is invariant, then \( \Lambda \frac{\mu_s(E_{\alpha})}{\mu(E_{\alpha})}=0 \) and \( \mu \sim \mu_{ac}+\mu_p \).
\end{theorem}
\begin{proof}[Proof]\par\nobreak
1. For the measure \( \gamma=\mu_s \) and an arbitrary \( m \in \mathbb{N} \), choose the set \( H \) and group elements \( g_1, ..., g_m \in G \) from Lemma \ref{lem:disjoint-translates}. The invariance of \( \varphi_{\mu} \) implies that
\[
m \varphi_{\mu}(H) = \sum\limits_{i=1}^m \varphi_{\mu}(g_iH) = \varphi_{\mu} \left( \bigsqcup\limits_{i=1}^m (g_iH) \right) \leq 1.
\]
Hence \( \varphi_{\mu}(H) \leq \frac{1}{m} \). Therefore,
\[
\Lambda \frac{\mu_s(E_{\alpha})}{\mu(E_{\alpha})} = \Lambda \frac{\mu_s(H \cap E_{\alpha})}{\mu(E_{\alpha})} \leq \Lambda \frac{\mu(H \cap E_{\alpha})}{\mu(E_{\alpha})} = \varphi_{\mu}(H) \leq \frac{1}{m}.
\]
Since \( m \) was arbitrary, we obtain
\[
\Lambda \frac{\| \mu \llcorner E_{\alpha} - (\mu_{ac} + \mu_p) \llcorner E_{\alpha} \|}{\mu(E_{\alpha})} = \Lambda \frac{\mu_s(E_{\alpha})}{\mu(E_{\alpha})} = 0.
\]\par
2. Since for every \( \alpha < \beta \) we have \( (\mu_{ac} + \mu_p)(E_{\alpha}) \leq  (\mu_{ac} + \mu_p)(E_{\beta}) \), the limit \( \lim\limits_{\mathfrak{F}} (\mu_{ac}+\mu_p)(E_{\alpha}) \) exists and equals \( \sup\limits_{\alpha \in I} (\mu_{ac}+\mu_p)(E_{\alpha}) \). Suppose that \( \lim\limits_{\mathfrak{F}} (\mu_{ac}+\mu_p)(E_{\alpha}) < +\infty \). Then
\begin{gather*}
1 = \Lambda \frac{\mu(E_{\alpha})}{\mu(E_{\alpha})} = \Lambda \frac{(\mu_{ac}+\mu_p)(E_{\alpha})}{\mu(E_{\alpha})} + \Lambda \frac{\mu_s(E_{\alpha})}{\mu(E_{\alpha})} = \lim\limits_{\mathfrak{F}} \frac{(\mu_{ac}+\mu_p)(E_{\alpha})}{\mu(E_{\alpha})}\\
= \frac{\lim\limits_{\mathfrak{F}} (\mu_{ac}+\mu_p)(E_{\alpha})}{\lim\limits_{\mathfrak{F}} \mu(E_{\alpha})} = 0,
\end{gather*}
which is impossible. Consequently, \( \lim\limits_{\mathfrak{F}} (\mu_{ac}+\mu_p)(E_{\alpha}) = +\infty \) and the measure \( \mu_{ac}+\mu_p \) is admissible.\par
Applying Lemma \ref{lem:sufficient-for-equivalence}, we obtain \( \mu \sim \mu_{ac} + \mu_p \).
\end{proof}
Fix a group \( \mathcal{G} \) equipped with the discrete topology.
\begin{definition}\label{def:invariant-mean}
\cite{Mitchell1965} A functional \( M \in B^*(\mathcal{G}) \) is called an invariant mean on \( \mathcal{G} \) if the following conditions hold:
\begin{enumerate}
\item \( \| M \| \leq 1 \)
\item \( M(1) = 1 \)
\item \( M(x \mapsto f(gx)) = M(f) \) for every function \( f \in B(\mathcal{G}) \) and every translation \( g \in \mathcal{G} \).
\end{enumerate}
\end{definition}
Conditions 1) and 2) immediately imply that \( M \) is positive.
\begin{definition}\label{def:left-thick}
\cite{Mitchell1965} A set \( A \subset \mathcal{G} \) is called left thick if for every finite set \( F \subset \mathcal{G} \) there exists \( g \in \mathcal{G} \) such that \( Fg = \{ fg: f \in F \} \) is contained in \( A \).
\end{definition}
If the group \( \mathcal{G} \) is abelian, a special case of \cite[Theorem 7]{Mitchell1965} states that for every left thick set \( C \subset \mathcal{G} \) there exists an invariant mean \( m_C \) on \( \mathcal{G} \) such that \( m_C(\chi_C)=1 \).\par
Let \( G = \langle \mathbb{R}, + \rangle \) have its usual topology, let \( I = \mathbb{N} \), \( \mathscr{E} = \{ [-N,N] \}_{N \in \mathbb{N}} \), and \( \mathcal{G} = \mathbb{R} / \mathbb{Z} \). We construct an admissible measure \(\mu\) with invariant \(\varphi_\mu\) such that no countably additive admissible measure is equivalent to \(\mu\) and \( \Lambda \frac{\mu_p(E_{\alpha})}{\mu(E_{\alpha})} \neq 0 \).
\begin{example}\label{ex:no-countably-additive-equivalent}\leavevmode\par\nobreak
1. Consider the canonical quotient map \( q \) from \( \mathbb{R} \) to \( \mathbb{R} / \mathbb{Z} \). Denote the open ball of radius \( \varepsilon > 0 \) centered at \( x \in \mathbb{R} / \mathbb{Z} \) by \( B(x, \varepsilon) \). To avoid a conflict of notation, for Haar measure on \( \mathbb{R} / \mathbb{Z} \) normalized to have total mass one we use \( \lambda_H \), and for Lebesgue measure on \( \mathbb{R} \) we use \( \lambda_L \). Enumerate the rational points of \( [0, 1) \) as \( \mathbb{Q} \cap [0, 1) = \{ r_i \}_{i=1}^{\infty} \). Set
\[
E = \bigcap\limits_{k=1}^{\infty} \bigcup\limits_{i=1}^{\infty} B (q(r_i), 2^{-k-i}).
\]
By continuity of Haar measure,
\begin{gather*}
\lambda_H(E) = \lim\limits_{k \to \infty} \lambda_H \left( \bigcup\limits_{i=1}^{\infty} B(q(r_i), 2^{-k-i}) \right) \leq \lim\limits_{k \to \infty} \sum\limits_{i=1}^{\infty} 2^{-k-i+1}\\
= \lim\limits_{k \to \infty} 2^{-k+1} \sum\limits_{i=1}^{\infty} 2^{-i} = \lim\limits_{k \to \infty} 2^{-k+1} = 0.
\end{gather*}
Thus \( \lambda_H(E)=0 \).\par
2. For an arbitrary finite set \( F \subset \mathbb{R} / \mathbb{Z} \), put
\[
L_F =  \{ l: F+l \subset E \} = \bigcap\limits_{f \in F} \{ l: f+l \in E \} =  \bigcap\limits_{f \in F} (E-f).
\]
As the intersection of at most countably many dense \( G_{\delta} \) sets, \( L_F \) is itself a dense \( G_{\delta} \) set \cite[p. 97]{Rudin1987}. In particular, \( L_F \neq \emptyset \). Hence \( E \) is a left thick set. Denote the invariant mean taking the value \( 1 \) on \( \chi_E \) by \( m_E \), and the measure on \( \mathscr{B}([0, 1)) \) assigning to each set \( A \) the value \( m_E(\chi_{q(A)}) \) by \( \mu_E \). Note that if \( B \in \mathscr{B}([0,1)) \) and \( h \in \mathbb{R} \) satisfy \( B+h \subset [0,1) \), then \( \mu_E(B+h) = \mu_E(B) \).\par
3. Define a measure \( \mu \) on \( \mathscr{B}(\mathbb{R}) \) by
\[
\mu(A) = \sum\limits_{n \in \mathbb{Z}} \mu_E((A-n) \cap [0, 1)).
\]
Fix an arbitrary set \( A \in \mathscr{B}(\mathbb{R}) \) and an arbitrary translation \( \Delta \in (0, 1) \). For every \( n \in \mathbb{Z} \), we have
\begin{align*}
&\mathrel{\phantom{=}}\mu_E((A+\Delta-n) \cap [0, 1))\\
&=\mu_E((A-n) \cap [-\Delta, 1-\Delta) + \Delta)\\
&=\mu_E((A-n) \cap [0, 1-\Delta) + \Delta) + \mu_E((A-n) \cap [-\Delta, 0) + \Delta)\\
&=\mu_E((A-n) \cap [0, 1-\Delta)) + \mu_E((A-n) \cap [-\Delta, 0)+1)\\
&=\mu_E((A-n) \cap [0, 1)) - \mu_E((A-n) \cap [1-\Delta, 1))\\
&\quad{}+\mu_E((A-(n-1)) \cap [1-\Delta, 1)).
\end{align*}
Put \( I=\mathbb{N} \), \( \Lambda = \lim\limits_{\mathbb{F}} \), and \( E_N = [-N, N] \) for \( N \in \mathbb{N} \). Then
\begin{align*}
&\mathrel{\phantom{=}}\varphi_{\mu}(A+\Delta)\\
&=\lim\limits_{N, \mathbb{F}} \frac{1}{2N} \sum\limits_{n=-N}^{N-1} (\mu_E((A-n) \cap [0, 1))\\
&\quad{}+\mu_E((A-(n-1)) \cap [1-\Delta, 1))-\mu_E((A-n) \cap [1-\Delta, 1)))\\
&=\lim\limits_{N, \mathbb{F}} \frac{1}{2N} \biggl( \sum\limits_{n=-N}^{N-1} \mu_E((A-n) \cap [0, 1))\\
&\quad{}+\mu_E((A+N+1) \cap [1-\Delta,1))-\mu_E((A-N+1) \cap [1-\Delta, 1)) \biggr).
\end{align*}
Since
\begin{gather*}
\left| \lim\limits_{N,\mathbb{F}} \frac{1}{2N} \left(\mu_E((A+N+1) \cap [1-\Delta,1))-\mu_E((A-N+1) \cap [1-\Delta, 1)) \right) \right|\\
\leq \lim\limits_{N,\mathbb{F}} \frac{1}{N} = 0,
\end{gather*}
\( \varphi_{\mu}(A+\Delta) = \varphi_{\mu}(A) \). Thus \( \varphi_{\mu} \) is invariant.\par
4. By the definition of \( \varphi_{\mu} \), we obtain
\begin{gather*}
\varphi_{\mu}(q^{-1}(E)) = \lim\limits_{N, \mathbb{F}} \frac{1}{2N} \sum\limits_{n=-N}^{N-1} \mu_E((q^{-1}(E)-n) \cap [0,1))\\
=\lim\limits_{N, \mathbb{F}} \frac{1}{2N} \sum\limits_{n=-N}^{N-1} \mu_E(q^{-1}(E) \cap [0,1)) = \mu_E(q^{-1}(E) \cap [0,1)) = m_E(\chi_E)=1.
\end{gather*}
By countable additivity and translation invariance of Lebesgue measure,
\begin{gather*}
\lambda_L(q^{-1}(E)) = \sum\limits_{n \in \mathbb{Z}} \lambda_L(q^{-1}(E) \cap [n-1,n)) = \sum\limits_{n \in \mathbb{Z}} \lambda_L(q^{-1}(E) \cap [0,1))\\
= \sum\limits_{n \in \mathbb{Z}} \lambda_H(E) = \sum\limits_{n \in \mathbb{Z}} 0 = 0.
\end{gather*}\par
5. Now suppose that there exists a countably additive admissible measure \( \gamma \) equivalent to \( \mu \). Then, by Theorem \ref{thm:remove-singular-part}, there exists an admissible measure \( \tau \) equivalent to \( \mu \) and absolutely continuous with respect to Lebesgue measure. In that case, \( \varphi_{\tau} = \varphi_{\mu} \) vanishes on every set of Lebesgue measure zero. This contradicts the equalities \( \varphi_{\mu}(q^{-1}(E)) = 1 \) and \( \lambda_L(q^{-1}(E))=0 \). Hence no countably additive admissible measure \( \gamma \) is equivalent to \( \mu \).\par
Using Theorem \ref{thm:remove-singular-part} and the equalities \( \lambda_L(q^{-1}(E))=0 \), \( \varphi_{\mu}(q^{-1}(E)) = 1 \) once again, we get
\begin{gather*}
\Lambda \frac{\mu_p(E_{\alpha})}{\mu(E_{\alpha})} = \Lambda \frac{(\mu_s+\mu_p)(E_{\alpha})}{\mu(E_{\alpha})} \geq \Lambda \frac{(\mu_s+\mu_p)(q^{-1}(E) \cap E_{\alpha})}{\mu(E_{\alpha})} = \Lambda \frac{\mu(q^{-1}(E) \cap E_{\alpha})}{\mu(E_{\alpha})}\\
= \varphi_{\mu}(q^{-1}(E)) = 1.
\end{gather*}
Consequently, \( \Lambda \frac{\mu_p(E_{\alpha})}{\mu(E_{\alpha})} \neq 0 \).
\end{example}

\section*{\centering Acknowledgments}
The author thanks Nikita V. Borisevich for investigating
a special case of Theorem \ref{thm:remove-singular-part}
for countably additive measures on the real line.

\section*{\centering AI Declaration}
ChatGPT (OpenAI) conducted the initial literature search and assisted with the retrieval of information from publications, the discussion of possible simplifications of Example \ref{ex:no-countably-additive-equivalent} and the proof of Theorem \ref{thm:SakM}, and the identification of possible errors and inaccuracies in mathematical proofs throughout the manuscript. ChatGPT also assisted with the English translation and language editing. The author is responsible for the content of the manuscript.\par

\end{document}